\documentclass{amsart}
\theoremstyle{definition}
\numberwithin{equation}{section}
\input xypic
\newcommand{\A}{\mathcal{A}}
\newcommand{\B}{\mathcal{B}}
\newcommand{\F}{\mathcal{F}}
\newcommand{\Hi}{\mathcal{H}}
\newcommand{\X}{\mathcal{X}}
\newcommand{\Y}{\mathcal{Y}}
\newcommand{\eL}{\mathcal{L}}

\newcommand{\N}{\mathbb{N}}
\newcommand{\Z}{\mathbb{Z}}
\newcommand{\R}{\mathbb{R}}
\newcommand{\C}{\mathbb{C}}
\newcommand{\lra}{\longrightarrow}

\begin{document}
\baselineskip=18pt

\title{Arveson's extension theorem in $*$-algebras}

\author{G. H. Esslamzadeh}
\address{Department  of Mathematics, Faculty of Sciences, Shiraz University, Shiraz 71454, Iran}
\email{esslamz@shirazu.ac.ir}
\author{L. Turowska}
\address{Department of  Mathematical Sciences, Chalmers University of Technology and the University of Gothenburg, SE-412 96 Goteborg, Sweden}
\email{turowska@chalmers.se}

\thanks{}

\begin{abstract}
Arveson's extension theorem asserts that $\B(\Hi)$ is an injective object in the category of operator systems. Calling every self adjoint unital subspace of a unital $*$-algebra, a quasi operator system, we show that Arveson's theorem remains valid in the much larger category of quasi operator systems. This shows that Arveson's theorem as a non commutative extension of Hahn-Banach theorem, is of purely algebraic nature.
\end{abstract}

\maketitle

\section{Preliminaries}
A  fundamental theorem of Arveson asserts that a completely positive linear map from a unital self adjoint subspace $\X$ of a $C^*$-algebra $\A$ into $\B(\Hi)$, where $\Hi$ is an arbitrary Hilbert space, can be extended to a completely positive map from $\A$ into $\B(\Hi)$. This theorem is now known as a non commutative generalization of the most widespreadly used result in functional analysis, the Hahn-Banach extension theorem. The algebraic nature of the Hahn-Banach's  theorem raises a natural question that whether this is the case for the Arveson's theorem as well. This question which is the main motivation of the present work, will be answered  positively in this paper. In the first glance one would expect that some modifications of the classic argument of Arveson or one of its alternates should suffice to prove the algebraic from of Arveson's theorem. However this does not work and a new approach is needed. The injective property of $\B(\Hi)$ has  been a motivation behind  applications and several other generalizations of Arveson's theorem. For instance see [4, 5, 11]. However  all of these generalizations  are focused on the range of a completely positive map in the sense that they replace $\B(\Hi)$ in the Arveson's theorem with an injective object  of a larger category of functional analytic objects; while in the present paper we expand the category of objects of domain, by replacing the domain with a bounded unital $*$-algebra. This selection was motivated by nice properties of the bounded part of a unital $*$-algebra which were revealed in [6, 7]. From this point of view the present paper can be considered as a sequel to [6].

This paper is organized as follows. Later in this section, we introduce our notations. In  Section 2 we  prepare the ground for our main result, by extending some classical results on $C^*$-algebras to unital $*$-algebras. Then we demonstrate necessity of a new approach for establishing a generalization of Arveson's theorem in the algebraic context. In section 3 which is the main part of this paper, we show that $\B(\Hi)$ is injective in the category of $*$-algebras with an order unit. One of the main tools in our argument is Archimedeanization, a technique of [16], which has been adapted to the general case of $*$-vector spaces with not necessarily proper positive cones. We conclude the paper by showing the universal property of Archimedeanization.

Before proceeding further, let us set up our terminology.

Throughout all vector spaces are over the complex field and by $*$-vector space we mean a vector space $\X$ with an additive involution whose set of self adjoint elements is denoted with $\X_{sa}$. Moreover $M_n(\X)$ denotes the set of all $n\times n$ matrices over $\X$ with usual matricial operations and involution.
A {\it cone} in $\X$ is a nonempty subset $C\subseteq \X$ which is closed under addition and multiplication by non negative reals.
If moreover $C$ satisfies the identity $C\bigcap(-C)=\{0\}$, then it is called a {\it proper cone}. $C$ is called a {\it full cone} if $\X_{sa}=C-C$.

A cone $C\subseteq\X_{sa}$ in a $*$-vector space $\X$ naturally induces a reflexive and transitive but not necessarily antisymmetric relation $\leq$  on $\X_{sa}$ by $x\leq y$ if $y-x\in C$, $(x,y\in \X_{sa})$ which we call  the {\it quasi order}  induced by $C$ and we call $\X$ a {\it quasi ordered $*$-vector space}. If moreover $C$ is proper, then the induced relation is a partial order on $\X_{sa}$ and $\X$ is called an ordered $*$-vector space. An element $1_\X\in C$ is called an {\it order unit} for $\X$ if for each $x\in \X_{sa}$, there exists a positive real number $r$ such that $x\leq r1_\X$. The cone $C$ is called {\it Archimedean} if for each $x_0\in C$, the following condition holds: If $x\in\X_{sa}$ and $x+rx_0\geq 0$ for every positive real number $r$, then $x\geq 0$.
If $\X$ has an order unit $1_\X$, one needs only to consider the element $x_0=1_\X$ for proving that $C$ is Archimedean. In this case $1_\X$ is called an Archimedean order unit.

 We say that $(\X, \{ C_n\}_{n=1}^\infty)$ is a {\it matrix quasi ordered $*$-vector space} if $\X$ is a $*$-vector space in which the following conditions hold:

(i) $C_n$ is a cone in $M_n(\X)$ for each $n$;

(ii) For each $m,n\in \N$ and each $\Lambda\in M_{n,m}(\C)$ we have $\Lambda^* C_n\Lambda\subseteq C_m$ (Compatibility condition).

Moreover if $\X$ contains an element $1_\X$ such that $1_\X\otimes I_n$ is an order unit for $M_n(\X)$, $(n\in\N)$, then $1_\X$ is called a matrix order unit.

Throughout $\A$ is a complex unital $*$-algebra.    An element $a$ of $\A$ is called {\it positive} if $a=\sum_{k=1}^n a_k^*a_k$ for some  $n\in\N$ and  some $a_1,\ldots,a_n\in\A$. The set of positive  elements of $\A$ is a cone in the real vector space $\A_{sa}$ which we denote by $\A^+$ and we call it the {\it algebraic positive cone} of $\A$.
An element $x$ of  $\A$ is said to be {\it bounded} if there exists a positive real number $k$ such that $x^*x\leq k1$. The set of all such elements is a $*$-subalgebra of $\A$ [2, \S 54, Proposition 1] which is denoted by $\A_0$ and  is called the {\it bounded subalgebra} of $\A$. Also we call the set of bounded elements $x$ with $x^*x\leq 1$  the {\it unit ball} of $\A$ and we denote it by $ball(\A)$. The  {\it reducing ideal} $A_R$  of $\A$ is the intersection of the kernels of all irreducible $*$-representations of $\A$ by bounded operators on some Hilbert space.

 If $\X$ is a subspace of $\A$, then by $\X^+$ we mean $\X\cap [\X]^+$ where $[\X]$ is the unital self adjoint subalgebra of $\A$ generated by $\X$; while by $\X_0$ we mean $\X\cap\A_0$ and by $ball(\X)$ we mean $\X\cap ball(\A)$. We denote by $\X^*$ the space $\{x^*:x\in\X\}$.

It was shown in [6, Theorem 2.8] that the bounded $*$-subalgebra $\A_0$  of $\A$ has a $C^*$-seminorm defined by
$$
\|a\|_{\A} = \inf\{k\in\R^+ : a^*a\leq k^21_{\A}\}, \ \ \ (a\in\A_0).
$$
Furthermore, if we assume that $1_{\A}$ is Archimedean for $\A_0$
and the cone $\A_0^+$ is proper, then $\|.\|_\A$ is a $C^*$-norm on $\A_0$.

Let $\X$ and $\Y$ be  vector spaces, $\eL(\X,\Y)$ be the space of linear maps from $\X$ into $\Y$ and $\phi\in\eL(\X,\Y)$. As usual for each $n\in\N $ the $n$-th amplification of $\phi$ is the map $\phi_n:M_n(\X)\longrightarrow M_n(\Y)$ given by $\phi_n([x_{ij}])=[\phi(x_{ij})]$. Let $\X$ and $\Y$ be quasi ordered $*$-vector spaces with cones $C$ and $D$ respectively. A linear map $\phi:\X\lra\Y$ is called {\it positive} [resp. {\it self adjoint}] if $\phi(x)\in D$ for all $x\in C$ [resp. $\phi(x^*)=\phi(x)^*$ for all $x\in\X$]. If $\phi_n$ is positive for all $n\in\N $, we say that $\phi$ is {\it completely positive}. We denote the space of all completely positive linear maps from $\X$ into $\Y$ by $CP(\X ,\Y)$.

Now we are in the position to define the key concepts of this article. These notions were introduced and studied in [6, 7].

{\bf Definition.}  Let $\A$ be a unital $*$-algebra. We call any subspace [resp.  unital self adjoint  subspace/ subalgebra] $\X$  of $\A$, a {\it quasi operator space} [resp. {\it quasi operator system/ quasi operator algebra}] in $\A$ or briefly a  quasi operator space [resp. quasi operator system/ quasi operator algebra] when the algebra $\A$ is clear from the context.

We will show in Lemma 3.1 that every bounded quasi operator space is a matrix quasi ordered $*$-vector space.

{\bf Definition.} Let  $\X$ and $\Y$ be  quasi operator spaces [resp. quasi operator systems] in the unital $*$-algebras  $\A$ and $\B$   respectively. A linear map $\phi:\X\lra\Y$ is called  {\it contractive} if $\phi(x)\in ball(\Y)$ for all $x\in ball(\X)$. Also we call $\phi$ {\it bounded} if there is an $r\in\R^+$ such that $\phi(x)^*\phi(x)\leq r1_{\B}$ for all $x\in ball(\X)$. If $\phi$ is a bounded linear map, then we define the {\it algebraic bound}  of $\phi$ as follows:
$$
|||\phi||| = \inf\{r\in\R^+\  :\  \phi(x)^*\phi(x)\leq r^21_{\B}\    \text{for all} \  x\in ball(\X)\}.
$$
 Furthermore, we say that $\phi$ is {\it completely bounded} if $\phi_n$ is bounded for all $n\in\N $ and $\sup\{|||\phi_n||| : n\in\N \}$ is finite which we call  the {\it algebraic completely bound} of $\phi$ and we denote this quantity by $|||\phi|||_{cb}$.

We use the notation $|||.|||$ to avoid confusion with usual operator norm which we denote with $\| .\|$. Let $B_{alg}(\X,\Y)$ be the set of all algebraically bounded linear maps from $\X$ into $\Y$. Using the fact that the bounded subalgebra of any complex unital $*$-algebra is itself a $*$-algebra [2, page 243], we see that  $B_{alg}(\X,\Y)$ is a complex vector space and the algebraic bound  $|||.|||$ defined above is a seminorm on it. Indeed when $\X=\Y$, by [2, Proposition 1, page 243]  $B_{alg}(\X,\Y)$ which we denote simply with    $B_{alg}(\X)$  is a unital algebra.  The set of all completely bounded maps from
$\X$ into $\Y$ is a vector space equipped with $|||.|||_{cb}$ which we denote with  $CB_{alg}(\X,\Y)$. Also we denote the set of {\it states} of $\X$, that is, positive linear functionals $f$  on $\X$ with $|||f|||=1$ by $S(\X)$.

{\bf Remarks 1.1.} (i) Let $\X$ [resp. $\Y$] be a quasi operator space in the unital $*$-algebra $\A$ [resp. $\B$] and $T\in B_{alg}(\X,\Y)$. Observe that the algebraic bound behaves like the operator norm and the following identity holds.
$$
|||T|||=sup \{ \|Tx\|_\Y \ :\  x\in ball(\X)\}.
$$
Moreover $ \|Tx\|_\Y \leq |||T||| \|x\|_\X$ for all $x\in \X_0$.

(ii) With the above notations when $\X$ and $\Y$ are operator spaces,  our definition of boundedness coincides with the usual notion of boundedness and algebraic bound of $\phi$ equals the operator norm of $\phi$.

A positive linear functional $f$ on $\A$  is called {\it admissible} if for every $x\in\A$ the linear map $\pi(x)\ :\ \A/N\lra\A/N$, $\pi(x)(y+N)=xy+N$ is bounded where  $N=\{ x\in\A\ : f(x^*x)=0\}$. Throughout by  $*$-representation we mean a $*$-homorphism into $\B(\Hi)$ for some Hilbert space $\Hi$. If for every $x\in\A$,
$$
sup\{ \|\pi(x)\|\ : \ \pi\  {\rm  is\  a\ *-representation\ of\ \A \ on\ some\ Hilbert\  space}\} <\infty
$$
then we say that $\A$ is a {\it $G^*$-algebra} and the above supremum defines a $C^*$-seminorm on $\A$ which is called the {\it Gelfand-Naimark seminorm} and is denoted with $||.||_\gamma$. A $G^*$-algebra on which every positive linear functional is admissible, is called a {\it $BG^*$-algebra}. $\A$ is called a {\it $T^*$-algebra} if for every self adjoint element $a\in\A$, there is a $t\in\R$ satisfying  $t1_\A+a\in\A^+$. It was shown in [6, Theorem 2.14] that $\A_0$ is the largest $T^*$-subalgebra of $\A$.

\section{Extensions of completely positive maps into $M_n$}

This is a preparatory section for building an algebraic analog of Arveson's theorem. To do so,  we  construct a version of  Krein's theorem for  unital $*$-algebras in Theorem 2.1, and then an algebraic version of Arveson's theorem for maps into $M_n$.

{\bf Theorem 2.1.} {\it Let $\X$ be a quasi operator system in $\A$ and $f$ be a positive linear functional on $\X$. Then $f|_{\X_0}$ can be extended to a positive linear functional $\widetilde{f}$ on $\A_0$ such that $|||\widetilde f|||=|||f|_{\X_0}|||$.}

{\it Proof.} One can easily deduce from the proof of  [6, Theorem 3.8] that $|||f|||=f(1)$. If  $|||f|||=0$, then by Remark 1.1 (i) we have $f(x)=0$ for all $x\in\X_0$. So $f|_{\X_0}=0$ and the zero functional is the obvious extension of $f$ with desired properties. Thus without loss of generality we may assume that $|||f|||=1$. By Remark 1.1 (i), for every $x\in\X_0$ we have $|f(x)|\leq\| x\|_{\A}$ and hence by the Hahn-Banach Theorem there is a linear functional $\widetilde f$ on $\A_0$ such that $\widetilde f|_{\X_0}=f$ and $|\widetilde f(x)|\leq\| x\|_\A$ for all $x\in\A_0$. Now these conditions imply that $\widetilde f$ is bounded,  $|||\widetilde f|||\leq 1$ and $\widetilde f(1)=f(1)$.  But  $1=|||f|||\leq |||\widetilde f|||$. Thus $ |||\widetilde f|||=1=\widetilde f(1)$ and by [6, Theorem 3.8],  $\widetilde f\geq 0$.  Therefore $\widetilde f$ has the desired properties. \qed

In what follows, we adapt  notations of [14].
Let  $\X$ be a quasi operator space in $\A$, $\{ e_j\}$ be the canonical basis for $\C^n$ and $\{E_{ij}\ :\ i,j=1,...n\}$ be the standard matrix unit system in $M_n$. If $\phi:\X\lra M_n$ is a linear map, associate to $\phi$ a linear functional  $f_\phi$ on $M_n(\X)$ via
$$
f_\phi[x_{ij}]=\frac 1n\sum_{i,j=1}^n <\phi(x_{ij})e_j\  ,\ e_i>,\quad ([x_{ij}]\in M_n(\X)).
$$
Conversely if $f$ is a linear functional on $M_n(\X)$, then associate to $f$ a linear map  $\phi_f:\X\lra M_n$ by the formula
$$
\phi_f(x)= n\sum_{i,j=1}^n f(x\otimes E_{ij}) E_{ij}, \quad (x\in\X).
$$
As in the classic case it is easy to see that the maps $\phi\lra f_\phi$ from $\eL(\X, M_n)$ into  $\eL(M_n(\X),\C)$ and $f\lra \phi_f$ from $\eL(M_n(\X),\C)$ into  $\eL(\X, M_n)$ are mutual inverses. Moreover if $\X$ is unital, then in this correspondence unital maps correspond to unital functionals.

The following Lemma is the algebraic version of [14, Lemma 3.13] which can be proved with appropriate changes in the argument of the aforementioned.

{\bf Lemma 2.2.} {\it Positive elements of $M_n(\A)$ are precisely those $[x_{ij}]$ which have  a representation
$$
[x_{ij}]=\sum_{k=1}^m\sum_{l=1}^n[y^*_{kli}y_{klj}]_{i,j=1}^n
$$
for some $\{ y_{kli}\}\subseteq\A$, that is, $[x_{ij}]$ is a sum of $mn$ elements of the from $[y_i^*y_j]$ for some $\{ y_1,\dots , y_n\}\subseteq\A$.}

{\bf Lemma 2.3.} {\it For every quasi operator space  $\X$  in $\A$ we have  $M_n(\X_0)_0=M_n(\X_0)= M_n (\X)_0$.}

{\it Proof.} Suppose $[x_{ij}]\in M_n(\X_0)$. Then for every $i,j\leq n$ there is a $k_{ij}\in\R^+$ such that $ k_{ij}^2 1_\A-x_{ij}^* x_{ij}\geq 0$. So
 \begin{eqnarray}
k_{ij}^2 1_\A\otimes I_n&-&(x_{ij}\otimes E_{ij})^*(x_{ij}\otimes E_{ij})=
(k_{ij}^2 1_\A\otimes I_n)-(x_{ij}^*x_{ij}\otimes E_{jj})\nonumber\\
&=& diag (k_{ij}^2 1_\A,...,k_{ij}^2 1_\A,k_{ij}^2 1_\A-x_{ij}^*x_{ij},k_{ij}^2 1_\A,...,k_{ij}^2 1_\A)\nonumber
\end{eqnarray}
where the last term represents a diagonal matrix  whose j-th diagonal entry is $k_{ij}^2 1_\A-x_{ij}^*x_{ij}$ and   $k_{ij}^2 1_\A$ elsewhere, which is non-negative. Thus $x_{ij}\otimes E_{ij}\in M_n(\A)_0$ and hence  $[x_{ij}]$ is bounded, since by [2, Proposition 1, page 243] $M_n(\A)_0$ is a $*$-subalgebra of $M_n(\A)$. Therefore $M_n(\X_0)\subseteq M_n(\X)_0$.

Conversely let $[x_{ij}]\in M_n(\X)_0$ and $r\in\R^+$ be such that $r^2 1_\A\otimes I_n-[x_{ij}]^*[x_{ij}]\geq 0$. It is easy to see that the diagonal entries of a positive element of $M_n(\A)$ are all positive. So for arbitrary $i,j\leq n$ we have $r^2 1_\A-\sum_l x_{lj}^*x_{lj}\geq 0$ and hence $r^2 1_\A-x_{ij}^*x_{ij}\geq\sum_{l\neq i} x_{lj}^*x_{lj}\geq 0$. Thus $x_{ij}\in \X_0$ for all $i,j\leq n$. Therefore $M_n(\X_0)=M_n(\X)_0$ which involves the equality $M_n(\X_0)_0=M_n(\X_0)$ as an immediate consequence. \qed

{\bf Remark 2.4.} We can not replace $M_n$ in the identity $M_n\otimes\X_0= (M_n\otimes\X)_0$ in Lemma 2.3 with an arbitrary unital $*$-algebra or even a simple infinite dimensional $C^*$-algebra. Let $b(\A)$ be the bounded part  of a pro-$C^*$-algebra $\A$ in the sense of [16], which is the set
$$
\{ a\in\A :  \|a\|_\infty=sup\{ p(a) :  p\ \text{is a continuous}\ C^*-\text{seminorm on}\ \A\}<\infty\}.
$$
As it was observed in [6, Example 5.2] $b(\A)$ coincides with $\A_0$. Now if $\A$ is the algebra of compact operators on a separable infinite dimensional Hilbert space and $\B=C(\Z^+)$, then by
[17, Remark 3.6],
$$
 (\A\otimes\B)_0=b(\A\otimes\B)\neq \A\otimes b(\B)=\A\otimes\B_0 .
$$

The following theorem  is an algebraic analog of [14, Theorem 6.1] whose proof can be obtained by using Lemmas 2.2, 2.3 and making appropriate modifications in the argument of  [14, Theorem 6.1].

{\bf Theorem 2.5.} {\it Let $\X$ be a quasi operator system in $\A$ and $\phi\in \eL(\X, M_n)$. The following statements are equivalent.

\noindent (i) $\phi|_{\X_0}$ is completely positive.

(ii) $\phi|_{\X_0}$ is n-positive.

(iii) $f_\phi|_{M_n(\X_0)}$ is positive.}

The following theorem which is a consequence of Theorem 2.5, is an algebraic analog of [14, Theorem 6.2] and  can be proved with a similar argument.

{\bf Theorem 2.6.} {\it Let $\X$ be a quasi operator system in $\A$ and $\phi\ :\ \X\lra M_n$ be completely positive. Then there is a completely positive map   $\psi\ :\ \A_0\lra M_n$ which extends $\phi|_{\X_0}$.}

The following Lemma which is an analog of [14, Proposition 2.12], can be proved by adapting the argument of  the aforementioned and applying [6, Lemma 3.4].

{\bf Lemma 2.7.} {\it Let $\X$ be a unital quasi operator space in $\A$ and $\B$ be a unital $C^*$-algebra. If $\phi\ :\ \X\lra\B$ is a unital contraction, then the map $\tilde\phi\ :\ \X_0+\X_0^*\lra\B$, $\tilde\phi(x+y^*)=\phi(x)+\phi(y)^*$ is the unique positive extension of $\phi|_{\X_0}$ to $\X_0+\X_0^*$.}

The following theorem is an algebraic generalization of [14, Theorem 6.3].

{\bf Theorem 2.8.} {\it Let $\X$ be a unital quasi operator space in $\A$ and let $\phi\ :\ \X\lra M_n$ be a unital linear map. The following conditions are equivalent:

(i) $\phi|_{\X_0}$ is completely contractive.

(ii) $\phi|_{\X_0}$ is n-contractive.

(iii) $f_\phi|_{M_n(\X_0)}$ is contractive.

In particular if $\phi|_{\X_0}$ is n-contractive, then it can be extended to a completely positive map on $\A_0$.}

{\it Proof. } We need only to prove (iii) $\Longrightarrow$ (i). Let $\Y=\X_0+\X_0^*$.  Then $\Y$ is a quasi operator system and the unique positive extension of $f_\phi|_{M_n(\X_0)}$   to $M_n(\Y)$ obtained by Lemma 2.7 which we denote  by $g$, is  unital. Now by Theorem 2.5, $g$ is associated with a completely positive map $\psi$ on $\Y$ which clearly should satisfy the identity $\psi(x+y^*)=\phi(x)+\phi(y)^*$. Therefore $\psi$ is the natural extension of $\phi|_{\X_0}$ to $\Y$ and  hence
$\phi|_{\X_0}$ is completely contractive, by [6, Theorem 4.3]. \qed

We can prove a generalization of Arveson's extension theorem in the  Banach $*$-algebra context directly, by adapting the classical  proof presented in  [13]. We do this in the next theorem. Indeed all proofs of  Arveson's theorem are of analytic nature and it seems that due to the key role of weak* compactness, they are not  directly applicable to general $*$-algebras. See Remark 2.10.

{\bf Definition.} Let $\A$ be a unital Banach $*$-algebra. We call every (not necessarily closed) subspace [resp. unital self-adjoint subspace] of $\A$ a normed quasi operator space [resp. normed quasi operator system] in $\A$.

Suppose $\X$ and $\Y$ are normed quasi operator spaces in the unital Banach $*$-algebras $\A$ and $\B$ respectively. It is known that $B_{alg}(\X,\Y')$  and $\X\widehat\otimes\Y$ form a dual pair. See for instance [14, Lemma 7.1] (Note that in [14, Lemma 7.1]  completeness of  $\X$ and $\Y$ are not used in the proof and the same argument works for the general case). The weak* topology induced on $B_{alg}(\X,\Y')$ via identification with $\X\widehat\otimes\Y$ is called the bounded weak topology. Now if $\Hi$ is a Hilbert space, then by $CP_r(\X,\Hi)$ we mean the set
$$
CP_r(\X ,\Hi)=\{ T\in\B(\X,\B(\Hi))\ :\ T\ \text{is}\ \text{completely}\ \text{positive}\ \text{and}\ |||T|||\leq r\} .
$$
{\bf Theorem 2.9.} {\it Let $\X$ be a normed quasi operator system in the unital Banach $*$-algebra $\A$ and $\phi\ :\ \X\lra \B(\Hi)$ be completely positive. Then there is a completely positive map   $\psi\ :\ \A\lra  \B(\Hi)$ which extends $\phi$.}

{\it Proof.} Let $\F$ be a finite dimensional subspace of $\Hi$ and $\phi_\F\ :\ \X\lra\B(\F)$ be the compression of $\phi$ to $\F$, i.e., $\phi_\F(x)\xi=P_\F\phi(x)\xi$, $\xi\in\F$, where $P_\F$ is the projection onto $\F$.
As it was shown in [6] $\A=\A_0$ and consequently $\X=\X_0$. Moreover $\B(\F)$ is isomorphic to $M_n$ for some $n\in\N$ and hence Theorem 2.6 is applicable to $\phi_\F$. Thus there is a completely positive map $T_\F\ :\ \A\lra\B(\F)$ which extends $\phi_\F$.
For each $a\in\A$ define $\psi_\F(a)$ to be $0$ on $\F^\perp$ and be equal to $T_\F(a)$ on $\F$. Then by [6, Theorem 4.3]  $\psi_\F$ is completely bounded and $\psi_\F\in CP_r(\X ,\Hi)$ where $r=\|\phi(1_\A)\|$. Next, observe that [14, Theorem 7.4] is valid if we replace operator system with normed quasi operator system, with the same argument. Therefore $CP_r(\X ,\Hi)$ is compact in the bounded weak topology. Now the rest of proof is as in [14, Theorem 7.5]. \qed

{\bf Remark 2.10.} Suppose $\X$ is a  quasi operator space in $\A$ and $\Hi$ is an arbitrary Hilbert space. We can define the projective tensor seminorm on $\X_0\otimes\B(\Hi)_*$ where $\B(\Hi)_*$ is the predual of $\B(\Hi)$, exactly in the same way as it is defined on tensor product of normed spaces and then taking the completion to form  the projective tensor product $\X_0\widehat\otimes\B(\Hi)_*$.  As in the classic argument we equip  $CP_r(\X,\Hi)$ with the bounded weak topology which is the relative weak* topology inherited from $B_{alg}(\X,\B(\Hi)_*)=B_{alg}(\X_0,\B(\Hi)_*)=(\X_0\widehat\otimes\B(\Hi)_*)'$. One of the crucial steps in the proof Arveson's theorem and in the previous theorem is the compactness of $CP_r(\X,\Hi)$ in the bounded weak topology which is a consequence of Alaoglu's theorem. But Alaoglu's theorem valid in topological vector spaces, while $\X_0\widehat\otimes\B(\Hi)_*$ is not a topological vector space unless it is Hausdorff. Now it is a standard fact that the topology induced by a single seminorm, namely the projective tensor seminorm, is Hausdorff if and only if the given seminorm is actually a norm. This forces the seminorm of $\X_0$ to be a norm too; A case which was covered in Theorem 2.9. Therefore we can not  adapt Arveson's proof or its alternates directly to a purely algebraic setting and a new approach is needed to establish the algebraic form of Arveson's theorem. This is achived in the next section, by using Archimedeanization technique.

\section{Arveson's extension theorem in Bounded  $*$-algebras}

In this section we prove  Arveson's extension theorem  in the bounded part of a $*$-algebra, by using the Archimedeanization technique. But first we need to establish some facts about the bounded part of a $*$-algebra, which are done in the next three lemmas.

{\bf Lemma 3.1.} {\it The bounded part $\A_0$ of $\A$ , equipped with the quasi order induced by the algebraic cone,  is a matrix quasi ordered $*$-vector space. Moreover if $1_\A$ is an order unit for $\A$, then it is a matrix order unit, that is $1_\A\otimes I_n$ is an order unit for $M_n(\A)$ for all $n\in\N$.}

{\it Proof.} For the first statement we need only to check the compatibility condition. By Lemma 2.3, positive elements of $M_m(\A)$ are precisely finite sums of elements of the from $[y_i^*y_j]$ where $y_1,\dots ,y_m\in\A$ are arbitrary. So it is enough to check the compatibility condition for these later elements.  Suppose $\Lambda=[\lambda_{ij}]\in M_{m\times n}$ and  $\Lambda^* [y_i^*y_j]\Lambda=[x_{ij}]$.   Then
$$
x_{ij}=\sum_{k,l}\overline\lambda_{ki}y_k^*y_l\lambda_{lj}=\sum_k \overline\lambda_{ki}y_k^*\sum_l y_l\lambda_{lj}=z_i^*z_j
$$
where  $z_i=\sum_{k=1}^m\lambda_{ki}y_k$. So by Lemma 2.2 we have $[x_{ij}]\in M_n(\A)^+$.

Now suppose $1_\A$ is an order unit for $\A$. By [6,  Theorem 2.14] this is equivalent to the boundedness of $\A$, that is, $\A=\A_0$. Now by Lemma 2.3, $M_n(\A_0)=M_n(\A)_0$ for all $n\in\N$. Thus invoking  [6,  Theorem 2.14] once again, we see that $1_\A\otimes I_n$ is an order unit for $M_n(\A)$ for every $n\in\N$, or equivalently,  $1_\A$ is a matrix order unit.  \qed

Next lemma was proved implicitly in [6, Theorem 2.14].

 {\bf Lemma 3.2.} {\it On $\A_0$ the Gelfand-Naimark seminorm is equal to $\|.\|_\A$.}

The following lemma is a generalization of [16, Lemma 3.14] to matrix quasi ordered $*$-vector spaces. The general case can be proved with the same argument, since the proof of [16, Lemma 3.14] does not depend on properness of  positive cone. So we do not repeat the proof here.

{\bf Lemma 3.3.} {\it Let $(\X ,\{ U_n\}_{n=1}^\infty ,1_\X)$ be a matrix quasi ordered $*$-vector space with matrix order unit $1_\X$,
$$
N_n=\bigcap\{ ker f\ :\ f\in S(M_n(\X))\},\quad  (n\in \N)
$$
and $N=N_1$. Then for every $n\in\N$ we have $N_n=M_n(N)$. In particular  for every $[\lambda_{ij}]\in M_{m,n}$ we have $[\lambda_{ij}]^*M_m(N)[\lambda_{ij}]\subseteq M_n(N)$.}

{\bf  3.4. Archimedeanization of matrix quasi ordered $*$-vector spaces.} In this subsection we discuss the Archimedeanization of matrix quasi ordered $*$-vector spaces with a matrix order unit. This process was defined in [15, 3.2] for ordered $*$-vector spaces with an order unit,   and was extended to matrix ordered $*$-vector spaces with matrix order unit in [16, 3.1]. Here we follow the terminology of [16].

Let $(\X,\{ U_n\}_{n=1}^\infty, 1_\X)$ be a matrix quasi ordered $*$-vector space with a matrix order unit such that for each $n\in\N$ the cone $U_n$ is full. Let $N$ and $N_n$, $n\in\N$ be as in Lemma 3.3. We shall always  identify $M_n(\X/N)$ with $M_n(\X)/M_n(N)$ and we use the terms $M_n(N)$ and $N_n$ interchangably, in light of Lemma 3.3.
Now define
$$
D_n:=\{ [x_{ij}]\in M_n(\X)_{sa}\ :\ r1_\X\otimes I_n+ [x_{ij}]\in U_n\ \text{for all}\ r\in\R^+\}
$$
and
$$
(N_n)_\R :=\bigcap\{ ker f\ :\ f \ \text {is a real state on}\ M_n(\X)\} .
$$
Since  $U_n$ is a full cone, then the matricial version of Propositions 3.10 and 3.11 of [15] hold for $M_n(\X)$ with the same proof:

{\it If $f\ :\ M_n(\X)\lra\C$ is a $\C$-linear functional, then $f$ is positive if and only if}
$$
f(x)=g(Re(x))+ig(Im(x))\quad (x\in M_n(\X))
$$
{\it for some positive $\R$-linear functional $g$ on $M_n(\X)$. }

Applying this result, one can observe that the following identities hold.

(i) $N_n=(N_n)_\R+ i(N_n)_\R$;

(ii) $(N_n)_\R=D_n\cap-D_n$;

(iii) $(N_n)_\R=(N_n)_{sa}$.

 Since $N_n$ is a self adjoint subspace of $M_n(\X)$, then $M_n(\X)/N_n$ can be equipped  with a well defined  $*$-operation,
$$
([x_{ij}]+N_n)^*=[x_{ij}]^*+N_n,\quad ([x_{ij}]\in M_n(\X)).
$$
 One can see easily that
$$
M_n(\X/N)_{sa}=\{ [x_{ij}]+N_n\ :\  [x_{ij}]\in M_n(\X)_{sa}\} .
$$
If we define
\begin{eqnarray}
C_n:&=&\{ [x_{ij}]+N_n\in M_n(\X/N)_{sa}\ :\nonumber\\
 &&(r1_\X\otimes I_n+[x_{ij}])+N_n\in U_n+N_n\ \text{ for all}\ r\in \R^+\},\nonumber
\end{eqnarray}
then independence of proofs of [15, Theorem 2.35] and [16, Proposition 3.16] from properness of cones $U_n$, allows us to have those results with the same proofs. Besides, all of the cones $C_n$ are proper. Therefore,

{\it  $(\X/N,\{ C_n\}_{n=1}^\infty ,1_\X+N)$ is an Archimedean matrix ordered $*$-vector space with matrix order unit, which we denote with $(\X_{Arch}, \{ C_n\}, 1_\X+N)$ or simply $\X_{Arch}$ and we call it Archimedeanization of $(M_n(\X) ,\{ U_n\}_{n=1}^\infty, 1_\X)$.}

Moreover the following observation of [15, Remark 3.17] is true  for matrix quasi ordered $*$-vector spaces: Archimedeanization of a matrix quasi ordered $*$-vector space  $(\X,\{ U_n\}_{n=1}^\infty, 1_\X)$ is obtained by Archimedeanization of  $(\X, U_n, 1_\X\otimes I_n)$  at each matrix level.

 Finally, the following alternate expression of the cones $C_n$ is very useful and is used in the proof of the next theorem.

\begin{equation}\label{Cn}
C_n=\{ [x_{ij}]+N_n\ :\ [x_{ij}]\in D_n\} ,\ (n\in\N).
\end{equation}
Indeed  inclusion of the set  in the right side of the above identity in $C_n$ follows from the definition of $D_n$. To see the reverse inclusion, suppose $[x_{ij}]+N_n\in C_n$ and $s\in\R^+$ is given. Choose an $r\in (0,s)$.  By the definition of $C_n$ we have
$$
r1_\A\otimes I_n+[x_{ij}]+N_n\in U_n+N_n
$$
and hence there are $[y_{ij}]\in U_n$ and $[z_{ij}]\in N_n$ such that
$$
r1_\A\otimes I_n+[x_{ij}]= [y_{ij}]+[z_{ij}].
$$
Since $[x_{ij}]+N_n$ is self adjoint, then we may assume that $[x_{ij}]$ is self adjoint too and hence so is $[z_{ij}]$ in the above identity.  Thus $[z_{ij}]\in D_n$ as $(N_n)_{sa}\subseteq (N_n)_\R\subseteq D_n$. So $(s-r)1_\A\otimes I_n+[z_{ij}]\in U_n$. On the other hand $r1_\A\otimes I_n+[y_{ij}]\in U_n$ as both terms are in $U_n$. Therefore
$$
s1_\A\otimes I_n+[x_{ij}]=r1_\A\otimes I_n+[y_{ij}]+(s-r)1_\A\otimes I_n+[z_{ij}]\in U_n
$$
and hence $[x_{ij}]\in D_n$, as required.

Now we are prepared to state our main result.

{\it {\bf Theorem 3.5.} Let $\X$ be a quasi operator system in $\A$ and $\phi\ :\ \X\lra\B(\Hi)$ be a completely positive linear map, where $\Hi$ is an arbitrary Hilbert space. Then there is a completely positive linear map $\psi\ :\ \A_0\lra\B(\Hi)$ which extends $\phi|_{\X_0}$.}

{\it Proof.} Without loss of generality we may assume $\A=\A_0$. Consequently $\X=\X_0$ and $1_\A$ is a matrix order unit for $\X$ by Lemma 3.1. We break the proof into six steps. The commutative diagram at the end of proof illustrates these steps visually.

{\bf Step 1.} {\it Identifying the positive cones.}

 In this preliminay step, to make the proof self contained, we introduce the cones that are involved. Throughout for every quasi operator system $V$ in $\A$, by $V^+$ we mean the algebraic positive cone of $V$, that is, $V\cap[V]^+$. Let $n\in\N$. By [6, Lemma 3.4] the positive cones $M_n(\X)^+$ and $M_n(\A)^+$ are full cones. So the Archimedeanization process (3.4) is applicable to both of $\X$ and $\A$.  Let  $N_n$, $N$, $D_n$ and $C_n$ be defined as in 3.4. Replacing $\X$ with $\A$,  let $W_n$, $W$,  $E_n$ and $K_n$ be  analogs of  $N_n$, $N$, $D_n$ and $C_n$ respectively, that is,
$$
N_n:=\bigcap\{ ker f\ :\ f\in S(M_n(\X))\},
$$
$$
W_n:=\bigcap\{ ker f\ :\ f\in S(M_n(\A))\},
$$
$$
D_n:=\{ [x_{ij}]\in M_n(\X)_{sa}\ :\ r1_\A\otimes I_n+ [x_{ij}]\in M_n(\X)^+\ \text{for all}\ r\in\R^+\} ,
$$
$$
E_n:=\{ [x_{ij}]\in M_n(\A)_{sa}\ :\ r1_\A\otimes I_n+ [x_{ij}]\in M_n(\A)^+\ \text{for all}\ r\in\R^+\} ,
$$
\begin{eqnarray}
C_n:=&&\{ [x_{ij}]+N_n\in M_n(\X/N)_{sa}\ :\nonumber\\
 &&(r1_\A\otimes I_n+[x_{ij}])+N_n\in M_n(\X)^+ +N_n\ \text{ for all}\ r\in \R^+\}.\nonumber
\end{eqnarray}

\begin{eqnarray}
K_n:=&&\{ [x_{ij}]+W_n\in M_n(\A/W)_{sa}\ :\nonumber\\
 &&(r1_\A\otimes I_n+[x_{ij}])+W_n\in M_n(\A)^+ +W_n\ \text{ for all}\ r\in \R^+\}.\nonumber
\end{eqnarray}
By identity (1) of subsection 3.4 we have
\begin{equation}
C_n=\{ [x_{ij}]+N_n\ :\ [x_{ij}]\in D_n\},
\end{equation}
\begin{equation}
K_n=\{ [x_{ij}]+W_n\ :\ [x_{ij}]\in E_n\},\ (n\in\N).
\end{equation}
Throughout the proof we assume for every $n\in\N$ the spaces $M_n(\X/N)$ and $M_n(\A/W)$  are equipped with the Archimedeanization positive cones $C_n$ and $K_n$ respectively, unless stated otherwise.

{\bf Step 2.} {\it For all $n\in \N$ we have  $M_n(\X\cap W)=M_n(\X)\cap W_n=N_n$.}

{\it Proof of step 2.} Let $ [x_{ij}]\in M_n(\X)\cap W_n$ and $f\in S(M_n(\X))$. By Lemma 2.3, $M_n(\A_0)=M_n(\A)_0$. So we can apply Theorem 2.1 to extend $f$ into a state $\widetilde f$ on  $M_n(\A)$, as $\A=\A_0$. Thus $f( [x_{ij}])=\widetilde f( [x_{ij}])=0$ as  $ [x_{ij}]\in W_n$ and hence $ [x_{ij}]\in N_n$. Conversely let $ [x_{ij}]\in N_n$. Since for every $f\in S(M_n(\A))$ we have $f|_{M_n(\X)} \in S(M_n(\X))$,  then $f( [x_{ij}])=0$ and hence $ [x_{ij}]\in W_n$.   Thus $M_n(\X)\cap W_n=N_n,\ (n\in\N)$. In particular $\X\cap W=N$. Now this identity together with Lemma 3.3 imply that $M_n(\X\cap W)=M_n(\X)\cap W_n=N_n$ for all $n\in \N$.

{\bf Step 3.} {\it With the above notations, $W$ is equal to the reducing ideal $\A_R$ of $\A$. In particular $W_n$ is a $*$-ideal in $M_n(\A)$, $(n\in\N)$. Moreover $\| .\|_{\A/W}$ is a norm and if $\B$ is the completion of  $(\A/W, \| .\|_{\A/W})$, then $K_n\subseteq M_n(\B)^+$, $(n\in\N)$.}

{\it Proof of step 3.} By [6, Theorem 3.8] every positive linear functional $f$  on $\A$ is bounded and hence is continuous with respect to  $\|.\|_\A$. So by [13,  Proposition 9.4.12] $f$ is admissible. Consequently $f$ is representable by [13, Theorem 9.4.15], as $\A$ is unital.  This fact together with [13, Theorem 9.7.2] imply that
\begin{eqnarray}
\A_R&=&\cap\{ ker(f)\ :\ f\ \text{ is a representable positive linear functional on}\ \A\}\nonumber\\
&=&\cap\{ ker (f)\ :\  f\in S(\A)\}=W.\nonumber
\end{eqnarray}
 Thus $W$ is a $*$-ideal of $\A$ and $\A/W$ is a reduced unital $*$-algebra. On the other hand by Lemma 3.2 the Gelfand-Naimark seminorm of $\A$ is equal to $\| .\|_\A$. So applying [13, Theorem 10.1.3] we get
$$
W=\{ x\in\A\ :\ \|x\|_\A=0\}.
$$
 Moreover one can  easily deduce $(\A/W)_0=\A/W$ from the identity $\A=\A_0$. Thus by [6, Theorem 2.14] $\A/W$ is a reduced unital $T^*$-algebra and hence  $\| .\|_{\A/W}$ is a $C^*$-norm by [13, Corollary 10.1.8]. Now let $\B$ be the completion  of $(\A/W, \| .\|_{\A/W})$ and $[x_{ij}]+W_n\in K_n$. We just showed that  $W_n$ is a $*$-ideal in $M_n(\A)$, from which we get the identity $M_n(\A)^+ +W_n=M_n(\A/W)^+$. So for every $r\in\R^+$,
\begin{eqnarray}
r1_\B\otimes I_n+ ([x_{ij}]+W_n)&=&(r1_\A\otimes I_n+W_n)+ ([x_{ij}]+W_n)\nonumber\\
&=&(r1_\A\otimes I_n+ [x_{ij}])+W_n\in M_n(\A)^+ +W_n\nonumber\\
&\subseteq& M_n(\B)^+.
\end{eqnarray}
But the positive cone $M_n(\B)^+$ is Archimedean. Therefore $[x_{ij}]+W_n\in M_n(\B)^+$.

{\bf Step 4.} {\it The inclusion map of $\frac{\X+W}{W}$ into $\A/W$ induces a unital, completely positive, injective linear map $\iota\ :\ (\frac{\X+W}{W})_{Arch}\lra\A_{Arch}(=\A/W)$. Moreover if $\B$ is as in step 3 and $J_n$ is the Archimedeanization positive cone of $\frac{\X+W}{W}$, then
$$
M_n(\B)^+\cap M_n\left(\frac{\X+W}{W}\right)=J_n,\quad (n\in\N).
$$}

{\it Proof of step 4.} The first statement is a  consequence of the fact that $W$ is a $*$-ideal in $\A$. So by [6, Lemma 3.4] the positive cones $M_n(\frac{\X+W}{W})^+$ are full cones and hence  Archimedeanization process is applicable.  Let  $V_n$, $V$, $L_n$ and $J_n$ be analogs of  $N_n$, $N$, $D_n$ and $C_n$ respectively, when $\X$ is replaced with $\frac{\X+W}{W}$, that is,
$$
V_n:=\bigcap\left\{ ker f\ :\ f\in S\left(M_n\left(\frac{\X+W}{W}\right)\right)\right\},
$$
\begin{eqnarray}
L_n:=&&\left\{ [x_{ij}]+W_n\in M_n\left(\frac{\X+W}{W}\right)_{sa}\ :\nonumber\right.\\
&&\left. r1_\A\otimes I_n+ [x_{ij}]+W_n\in M_n\left(\frac{\X+W}{W}\right)^+\ \text{for all}\ r\in\R^+\right\},\nonumber
\end{eqnarray}

\begin{eqnarray}
J_n:=&&\left\{ ([x_{ij}]+W_n)+V_n\in M_n\left(\frac{(\frac{\X+W}{W})}{V}\right)_{sa}\ :\nonumber\right. \\
 &&\left. (r1_\A\otimes I_n+[x_{ij}]+W_n)+V_n\in M_n\left(\frac{\X+W}{W}\right)^+ +V_n\ \text{ for all}\ r\in \R^+\right\}.\nonumber
\end{eqnarray}
Let   $[x_{ij}]+W_n\in V_n$. Every $f\in S(M_n(\A))$ induces a state $\widetilde f\in S(M_n(\A/W))$, as $f(W_n)=0$. So $f[x_{ij}]=\widetilde f([x_{ij}]+W_n)=0$ and hence $[x_{ij}]\in W_n$. Thus $V_n=\{0\}$ and hence $J_n=L_n$. This fact together  with the relations
$$
M_n\left(\frac{\X+W}{W}\right)^+\subseteq M_n(\A/W)^+=M_n(\A)^++W_n
$$
imply that $J_n\subseteq K_n$. Therefore the inclusion map $\iota\ :\ \frac{\X+W}{W}\lra\A/W$ satisfies the required properties, as $V_n=\{0\}$, $(n\in\N)$.

As it was shown in step 3, $K_n\subseteq M_n(\B)^+$, which together with the relation $J_n\subseteq K_n$ leads to
$$
J_n\subseteq M_n(\B)^+\cap M_n\left(\frac{\X+W}{W}\right),\quad (n\in\N).
$$
Conversely, let $[x_{ij}]\in M_n(\B)^+\cap M_n\left(\frac{\X+W}{W}\right)$ and $f\in S\left(M_n\left(\frac{\X+W}{W}\right)_{Arch}\right)$. Since $M_n\left(\frac{\X+W}{W}\right)^+\subseteq J_n$, then $f\in S\left(M_n\left(\frac{\X+W}{W}\right)\right)$ as well. So by Theorem 2.1, $f$  can be extended to a state $g$ on $M_n(\A/W)$. But $(\A/W, \| .\|_{\A/W})$ is a unital operator algebra by step 3 and hence as it was pointed out in Remark 1.1, the algebraic bound $||| g|||$ is equal to the operator bound $\|g\|$. Therefore the continuous extension $\widetilde f$ of $g$ to $M_n(\B)$ satisfies the identity $\|\widetilde f\|=1=\widetilde f(1_\A\otimes I_n+W_n)$ and hence is a state on $M_n(\B)$.  So $f([x_{ij}])=\widetilde f([x_{ij}])\geq 0$. Now applying [15, Proposition 2.20] in the ordered real vector space $\left(M_n\left(\frac{\X+W}{W}\right)_{sa}, J_n\right)$, we conclude that $[x_{ij}]\in J_n$.

{\bf Step 5.} {\it $\X_{Arch}$ is completely order isomorphic to $(\frac{\X+ W}{W})_{Arch}$.}

{\it Proof of step 5.} By step 2 we can identify  $\X\cap W$ with $N$. On the other hand using notations of the previous step, we saw that $V_n=\{0\}$ and $J_n=L_n$. So we only need to prove that the algebraic isomorphism $T\ :\ \X/N\lra\frac{\X+ W}{W}$, $T(x+N)=x+W$ obtained from the second isomorphism theorem, is a complete order isomorphism, when $\X/N$ and $\frac{\X+ W}{W}$ are equipped with the Archimedeanization positive cones $\{C_n\}_{n=1}^\infty$ and $\{J_n\}_{n=1}^\infty$ respectively.

Suppose $ [x_{ij}]+N_n\in C_n$. By (\ref{Cn}) we may assume that $ [x_{ij}]\in D_n$. Let $r\in\R^+$. Then $r1_\A\otimes I_n+[x_{ij}]\in M_n(\X)^+$ and hence there is a finite set $\{ [z^1_{ij}],...,[z^m_{ij}]\}$ in the algebra $[M_n(\X)]$ generated by $M_n(\X)$ such that
$$
r1_\A\otimes I_n+[x_{ij}]=[z^1_{ij}][z^1_{ij}]^*+...+[z^m_{ij}][z^m_{ij}]^*.
$$
Keeping in mind that $W_n$ is a $*$-ideal of $M_n(\A/W)$ and elements of $[M_n(\X)]$ are finite sums of finite products of elements of $M_n(\X)$,  one can easily observe that $[z^l_{ij}]+W_n$ is in the algebra $\left[M_n\left(\frac{\X+ W}{W}\right)\right]$. Thus
$$
\left(r1_\A\otimes I_n+[x_{ij}]\right)+W_n=\sum_{l=1}^m([z^l_{ij}]+W_n)([z^l_{ij}]+W_n)^*\in\left[M_n\left(\frac{\X+ W}{W}\right)\right]^+
$$
and hence
$$
\left(r1_\A\otimes I_n+[x_{ij}]\right)+W_n\in M_n\left(\frac{\X+ W}{W}\right)^+.
$$
So $[x_{ij}]+W_n\in J_n$ as $r\in\R^+$ was arbitrary. Therefore $T$ is completely positive.

Conversely suppose $T_n([x_{ij}]+N_n)=[x_{ij}]+W_n\in J_n$ and $r\in\R^+$. Then
$$
\left(r1_\A\otimes I_n+[x_{ij}]\right)+W_n\in\left[M_n\left(\frac{\X+ W}{W}\right)\right]^+.
$$
Again applying the fact that  $W_n$ is a $*$-ideal of $M_n(\A/W)$ and recalling the algebraic structure of elements  of  $\left[M_n\left(\frac{\X+ W}{W}\right)\right]$, we can find a finite set $\{[y^1_{ij}],\dots ,[y^m_{ij}]\}$ in the algebra $[M_n(\X)]$  such that
\begin{eqnarray}
\left(r1_\A\otimes I_n+[x_{ij}]\right)+W_n=([y^1_{ij}]^*&+&W_n)([y^1_{ij}]+W_n)\nonumber\\
&+&\dots +([y^m_{ij}]^*+W_n)([y^m_{ij}]+W_n).\nonumber
\end{eqnarray}
Even more, there exists $[y_{ij}]\in [M_n(\X)]^+$ such that
$$
\left(r1_\A\otimes I_n+[x_{ij}]\right)+W_n=[y_{ij}]+W_n.
$$
On the other hand positivity of $\left(r1_\A\otimes I_n+[x_{ij}]\right)+W_n$ together with the identity ($\ref{Cn}$) preceding this theorem, allow us to assume that $[x_{ij}]$ is self adjoint. Now if we define $N'_n$, $C'_n$ and $D'_n$ to be analogs of $N_n$, $C_n$ and $D_n$ for $[M_n(\X)]$ respectively, then by step 2 we have $W_n\cap [M_n(\X)]=N'_n$ and hence by 3.4 we get
$$
\left(r1_\A\otimes I_n+[x_{ij}]-[y_{ij}]\right)\in (W_n\cap [M_n(\X)])_{sa}=(N'_n)_{sa}=(N'_n)_\R=D'_n\cap -D'_n.
$$
In particular for every $s\in\R^+$ we have
$$
s1_\A\otimes I_n+(r1_\A\otimes I_n+ [x_{ij}]-[y_{ij}])\in[M_n(\X)]^+
$$
and hence $(s+r)1_\A\otimes I_n+ [x_{ij}]\in [M_n(\X)]^+$.  This together with the fact that
$(s+r)1_\A\otimes I_n+ [x_{ij}]\in M_n(\X)$ imply that
$$
(s+r)1_\A\otimes I_n+ [x_{ij}]\in M_n(\X)\cap[M_n(\X)]^+=M_n(\X)^+.
$$
Since $r$ and $s$ were arbitrary positive real numbers, then  $[x_{ij}]\in D_n$ and hence $[x_{ij}]+N_n\in C_n$. Therefore $T^{-1}$ is completely positive.

{\bf Step 6.} {\it Construction of extension.}

{\it Proof of step 6.} Let $[x_{ij}]\in D_n$. Then for every $r\in\R^+$, $r1_\A\otimes I_n+[x_{ij}]\in M_n(\X)^+$ and hence
$$
r\phi(1_\A)\otimes I_n+\phi_n([x_{ij}])\in M_n(\B(\Hi))^+,\ (r\in\R^+).
$$
So by positivity of $\phi(1_\A)\otimes I_n$ and the Archimedean property of $M_n(\B(\Hi))^+$ we have $\phi_n([x_{ij}])\in M_n(\B(\Hi))^+$.

Now applying the identities $(N_n)_\R=D_n\cap-D_n$ and $N_n=(N_n)_\R+ i(N_n)_\R$, we see that $\phi_n(N_n)=0$. So the map
$$
\theta_n\ :\ M_n(\X/N)\lra M_n(\B(\Hi)),\ \theta_n([x_{ij}]+N_n)=\phi_n([x_{ij}])
$$
is well defined. Moreover if $[x_{ij}]+N_n\in C_n$, then by the identity (1) of subsection 3.4 we may assume that $[x_{ij}]\in D_n$. Thus
 $$
\theta_n([x_{ij}]+N_n)=\phi_n([x_{ij}])\in M_n(\B(\Hi))^+.
$$
On the other hand all of these maps $\theta_n$ are compatible, that is, $\theta_n=(\theta_1)_n$, as this is the case for $\phi_n$s and hence we have a natural completely positive map
$$
\theta\ :\ \X_{Arch}\lra \B(\Hi),\ \theta(x+N)=\phi(x),\ (x\in\X).
$$
Now we can identify $\X_{Arch}$ with the subspace  $(\frac{\X+ W}{W})_{Arch}$ of $\A_{Arch}$ via the map $T$ of step 5 and getting a completely positive map
$$
\widetilde\phi=\theta\circ T^{-1}\ :\  (\frac{\X+ W}{W})_{Arch}\lra\B(\Hi).
$$
By step 4 the inclusion map $\iota\ :\ (\frac{\X+W}{W})_{Arch}\lra\A_{Arch}$ is completely positive. On the other hand we showed in step 3 that $K_n\subseteq M_n(\B)^+$, $(n\in\N)$ where $\B$ is the completion of  $(\A/W, \| .\|_{\A/W})$. So the inclusion map  of $\A_{Arch}$ into $\B$ is completely positive, whose composition with $\iota$ provides a completely positive linear map $\widetilde\iota\ :\ (\frac{\X+W}{W})_{Arch}\lra\B$. Therefore we can view $(\frac{\X+W}{W})_{Arch}$ as an operator system in $\B$. Now the identity
$M_n(\B)^+\cap M_n\left(\frac{\X+W}{W}\right)=J_n$ of step 4 implies that $\widetilde\phi$ is completely positive in the $C^*$ sense and we can apply the classical Arveson's extension theorem to obtain a completely positive map $\widetilde\psi\ :\ \B\lra\B(\Hi)$ which extends $\widetilde\phi$.

Since $W$ is a $*$ ideal in $\A$, then the quotient map $\pi\ :\ \A\lra\A/W$ is completely positive and hence $\psi:=\widetilde\psi\circ\pi$  is a completely positive map such that for every $x\in\X$,
$$
\psi(x)=\widetilde\psi(x+W)=\widetilde\phi(x+W)=\theta\circ T^{-1}(x+W)=\theta(x+N)=\phi(x) .
$$
Therefore $\psi$ is an extension of $\phi$, as required.

The following commutative diagram in which all maps are completely positive, is a visual summary of the present proof.
\newpage
$$
\diagram\ \X\ \xto[r]^\phi\xto[d]^P&\ \ \B(\Hi)\ \ &\ \A\ \xto[l]_\psi\xto[d]^\pi\\
\ \ \ \ \ \ \X_{Arch}\ \xto[ur]^\theta\xto[r]^T & \ \ (\frac{\X+ W}{W})_{Arch}\xto[u]^{\widetilde\phi}\xto[r]^{\widetilde\iota}&\ \B\ \xto[ul]_{\widetilde\psi}\enddiagram
$$
\qed

Archimedeaniztion of matricially ordered $*$-vector spaces with matrix order unit are universal objects in the category of matrix ordered $*$-vector spaces with matrix order unit; A property which was proved in [16, Theorem 3.18]. In the next theorem we prove a similar property of quasi operator systems whose proof, although different,  was inspired by  [16, Theorem 3.18].

{\bf Theorem 3.6.} {\it Suppose $\X$ is a bounded quasi operator system in $\A$ and $\X_{\rm Arch}$ is the Archimedeanization of $\X$. Then the quotient map $P\ :\ \X\lra\X_{\rm Arch}$ is a unital completely positive map. Moreover if $\Y$ is a quasi operator system in $\B$ where $\B_0^+$ is a proper Archimedean cone and $\phi\ :\ \X\lra\Y$ is a unital completely positive map, then there exists a unique completely positive linear map $\widetilde\phi\ :\ \X_{\rm Arch}\lra\Y$ such that $\phi=\widetilde\phi\circ P$, that is, the following diagram is commutative.}
$$
\diagram\  \X \xto[r]^\phi\xto[d]^P &\Y\\ \ \ \ \ \X_{Arch}\xto[ur]^{\widetilde\phi} & \enddiagram
$$

{\it Proof.}  One can see that $P$ is completely positive by the definition of Archimedeanization. Let $\Y$ be a quasi operator system in $\B$ and $\phi\ :\ \X\lra\Y$ be a unital completely positive map. By [6, Theorem 4.3] $\phi$ is completely bounded. On the other hand by steps 1 and 2 in the proof of Theorem 3.5, $N\subseteq W=\A_R$. So by Lemma 3.2, $\|x\|_\A=0$ for all $x\in N$. Now Remark  1.1 shows that $\|\phi(x)\|_\B=0$ for all $x\in N$. But by [6, Theorem 2.8] $\|.\|_\B$ is a norm. So $\phi(N)=0$  and hence the map $\widetilde\phi\ :\ \X_{\rm Arch}\lra\Y$, $\widetilde\phi(x+N)=\phi(x)$, is well defined. Let $n\in\N$ and $[x_{ij}]+N_n\in C_n$. Then
$$
(r1_\A\otimes I_n+ [x_{ij}])+N_n\in D_n+N_n
$$
for all $r\in\R^+$ and hence
$$
r1_\B\otimes I_n+ \widetilde\phi_n([x_{ij}])\in M_n(\Y)^+=M_n(\Y)\cap[M_n(\Y)]^+\subseteq M_n(\B)^+
$$
for all $r\in\R^+$. Now the Archimedean property of $1_\B\otimes I_n$ implies that
$\widetilde\phi_n([x_{ij}])\in M_n(\B)^+$. Therefore $\widetilde\phi$ is completely positive. If $\psi\ :\ \X_{\rm Arch}\lra\Y$ is another completely positive linear map which satisfies the identity $\phi=\psi\circ P$, then for every $x\in\X$ we have
$$
\psi(x+N)=(\psi\circ P)(x)=\phi(x)=(\widetilde\phi\circ P)(x)=\widetilde\phi(x+N) .
$$
Therefore $\widetilde\phi$ is unique. \qed

{\bf Acknowledgements.} The authors would like to express their sincere thanks to Professor I. G. Todorov for his valuable comments. This research was done when the first author was visiting Department of Mathematical Sciences of Chalmers University of Technology and the University of Gothenburg. The first  author wishes to express his sincere thanks to the Department and Professor M. Asadzadeh for their hospitality during this visit.



\begin{thebibliography}{30}
\bibitem{1} {R. Alizadeh and  G. H. Esslamzadeh}, {\it On the existence of unitarily invariant norm under some conditions}, Lin.  Multilin. Alg.  58 (3)  (2010), 367-375
\bibitem{2} {S. K. Berberian}, {\it Baer $*$-rings}, Springer Verlag, 1972.
\bibitem{3} {S. K. Berberian}, {\it Baer rings and Baer $*$-rings}, Springer Verlag, 2003.
\bibitem{4} {M. D. Choi and E. G. Effros}, {\it Injectivity and operator spaces},  J. Funct. Anal. 24 (1977), 156-209.
\bibitem{5} {A. Dosiev}, {\it Local operator spaces, unbounded operators and multinormed $C^*$-algebras},  J. Funct Anal. 255 (2008), 1724-1760.
\bibitem{6} {G. H. Esslamzadeh and F. Taleghani}, {Structure of quasi operator systems}, Lin. Alg. Appl. 438 (2013), 1372-1392.
\bibitem{7} {G. H. Esslamzadeh and F. Taleghani}, {A characterization of bounded quasi operator spaces}, Submitted.
\bibitem{8} {E. G. Effros and Z. J. Ruan}, {\it On the abstract characterization of operator spaces}, Proc. Amer. Math. Soc. 119 (1993), 579-584.
\bibitem{9} {E. G. Effros and Z. J. Ruan}, {\it Operator Spaces}, Clarendon Press, 2000.

\bibitem{10} {R. V. Kadison and J. R. Ringrose}, {\it Fundamentals of the theory of operator algebras}, Academic Press, 1983.
\bibitem{11} {Karn}, {\it A p-theory of ordered normed spaces}, Positivity 14 (2010), 441-458.
\bibitem{12} {G. J. Murphy}, {\it $C^*$-algebras and operator  theory}, Academic Press, 1991.
\bibitem{13} {T. W. Palmer}, {\it Banach algebras and the general theory of $*$-algebras II }, Cambridge University Press, 2001.
\bibitem{14} {V. I. Paulsen}, {\it Completely bounded maps and operator algebras}, Cambridge University Press, 2002.
\bibitem{15} {V. I. Paulsen and M. Tomforde}, {\it Vector spaces with an order unit}, Indiana Univ. Math. J. 58 (2009), no.3, 1319-1359.
\bibitem{16} {V. I. Paulsen, I. G. Todorov, and M. Tomforde}, {\it Operator system structures on ordered spaces}, Proc. London Math. Soc. 102(1) (2011),  25-49.
\bibitem{17} {N. C. Phillips}, {\it Inverse limits of $C^*$-algebras}, J. Oper. Theory 19 (1989),  159-195.

\bibitem{18} {Z. J. Ruan}, {\it Subspaces of $C^*$-algebras}, J. Funct. Anal. 76 (1988), 217-230.
\bibitem{19} {W. J. Schreiner}, {\it Matrix regular operator spaces}, J. Funct. Anal. 152 (1998), 136-175.
\bibitem{20} {R. R. Smith}, {\it Completely bounded maps between $C^*$-algebras}, J. London Math. Soc. 27 (1983), 157-166.
\end{thebibliography}
\end{document}